\title{Ternary cyclotomic polynomials having a large coefficient}
\author{Yves Gallot and Pieter Moree}
\documentclass[12pt]{article}
\usepackage{amssymb, latexsym, amsfonts}
\def\@ptsize{2}
\newtheorem{Thm}{Theorem}
\newtheorem{Con}{Conjecture}
\newtheorem{Lem}{Lemma}

\newtheorem{Prop}{Proposition}
\newcommand{\qed}{\hfill $\Box$}

\begin{document}
\date{}
\maketitle
{\def\thefootnote{}
}

\begin{abstract}
\noindent Let $\Phi_n(x)$ denote the $n$th cyclotomic polynomial. 
In 1968 Sister Marion Beiter conjectured that
$a_n(k)$, the coefficient of $x^k$ in $\Phi_n(x)$, satisfies
$|a_n(k)|\le (p+1)/2$ in case
$n=pqr$ with $p<q<r$ primes (in this case $\Phi_n(x)$ is
said to be ternary). Since then several results towards
establishing her conjecture have been proved (for
example $|a_n(k)|\le 3p/4$). Here we show that,
nevertheless, Beiter's conjecture is false for every $p\ge 11$. We also prove that
given any $\epsilon>0$ there exist infinitely many triples $(p_j,q_j,r_j)$
with $p_1<p_2<\ldots $ consecutive primes such that $|a_{p_jq_jr_j}(n_j)|>(2/3-\epsilon)p_j$
for $j\ge 1$.
\end{abstract}
\section{Introduction}
The $n$th cyclotomic polynomial $\Phi_n(x)$ is defined by
$$\Phi_n(x)=\prod_{j=1\atop (j,n)=1}^n(x-\zeta_n^j)=
\sum_{k=0}^{\varphi(n)}a_n(k)x^k,$$
where $\varphi$ is Euler's totient function
and $\zeta_n$ a primitive $n$th root of unity. For the $k$ not in the
range $[0,\varphi(n)]$, we put $a_n(k)=0$.
The coefficients $a_n(k)$
are known to be integers. The study of the $a_n(k)$ began with the
startling observation that for small $n$ we have $|a_n(k)|\le 1$ (it
thus seems, as D. Lehmer \cite{DL} worded it, that the primitive roots
of unity conspire to achieve this smallness). The first
counter-example to $|a_n(k)|\le 1$ occurs for $n=105$: $a_{105}(7)=-2$. Note
that $105$ is the smallest odd integer having three prime factors. If $\omega_1(n)$
denotes the number of odd prime factors of $n$, then it is well-known that
if $\omega_1(n)\le 2$, then $\Phi_n(x)$ is {\it flat}, that is all its
coefficients satisfy $|a_n(k)|\le 1$. Thus $n=105$ is the first candidate
integer for $\Phi_n(x)$ to be non-flat. We see that with respect to the 
smallness of the coefficients the first non-trivial case arises when
$\omega_1(n)=3$. In this case some authors say that $\Phi_n(x)$ is
{\it ternary}. Then we write $n=pqr$ with $2<p<q<r$.\\
\indent We define 
the
{\it height} of $\Phi_n(x)$ to be
$\max\{|a_n(k)|:0\le k\le \varphi(n)\}$ and
denote it by $A(n)$.
In 1968 Sister Marion Beiter \cite{Beiter-1} put forward the 
following conjecture (which she repeated in 1971 \cite{Beiter}). 
\begin{Con} {\rm (Sister Marion Beiter, 1968)}.
If $2<p<q<r$ are primes, then $A(pqr)\le {p+1\over 2}$.
\end{Con}
Note that $A(2qr)=1$.
In case either $q$ or $r\equiv \pm 1({\rm mod~}p)$ Beiter proved her conjecture.
This result was extended by Bachman \cite{B1}.
\begin{Thm} 
\label{flup}
{\rm (Gennady Bachman, 2003).}
If either $q$ or $r$ is 
congruent to $\pm 1$ or $\pm 2$ modulo $p$, then $A(pqr)\le (p+1)/2$. If $q$ or $r$ is 
congruent to
$(p\pm 1)/2$ modulo $p$, then $A(pqr)\le (p+3)/2$.
\end{Thm}
\indent In a further paper Beiter \cite{Beiter} points out that her conjecture is
true for $p\le 5$ and shows that $A(pqr)\le p-\lfloor p/4 \rfloor$, thus
improving on a result from Bang \cite{Bang} proved in 1895, to the
effect that 
$A(pqr)\le p-1$. Bloom \cite{Bloom} independently showed that $A(5qr)\le 3$ (and hence
the truth of Beiter's conjecture for $p=5$).
The best known general upper bound to date is due to Bachman \cite{B1},
who proved that $A(pqr)\le p-\lceil p/4 \rceil$. In the same paper Bachman
showed:
\begin{Thm} {\rm (Gennady Bachman, 2003).}
\label{flup2}
Let $q^*$ and $r^*$, $0<q^*,r^*<p$ be the
inverses of $q$ and $r$ modulo $p$ respectively. Set
$a=\min(q^*,r^*,p-q^*,p-r^*)$.
Then $A(pqr)\le \min({p-1\over 2}+a,p-a)$.
\end{Thm}
H. M\"oller \cite{HM2} indicated for every prime $p>3$ a cyclotomic polynomial $\Phi_{pqr}(x)$
having a coefficient equal to $(p+1)/2$. This shows that
Beiter's conjecture is best possible, if true. In particular M\"oller proved:
\begin{Thm} {\rm (Herbert M\"oller, 1971)}. Let $3<p<q<r$ be prime
numbers satisfying $q\equiv 2({\rm mod~}p)$ and $r=(mpq-1)/2$ for some
integer $m$. Then $$a_{pqr}({1\over 2}{(p-1)(qr+1)})={p+1\over 2}.$$
\end{Thm}
Earlier Emma Lehmer \cite{Emma} had shown that for $q$ and $r$ as in the latter Theorem
we have $a_{pqr}({1\over 2}{(p-3)(qr+1)})=(p-1)/2$. On combining M\"oller's result
with Theorem \ref{flup} we infer that for his choice of $p,q$ and $r$ we have
$A(pqr)=(p+1)/2$.\\
\indent Let $M(p)$ be the maximum of the heights of the ternary cyclotomic polynomials,
where $p$ is the smallest prime factor of $n$. The case $p=3$
was investigated in detail by Beiter \cite{Beiter3}, who found that $M(3)=2$.
Beiter's conjecture in combination
with $M(3)=2$ and M\"oller's result leads to the following conjecture.
\begin{Con}
For $p>2$ we have $M(p)={p+1\over 2}$.
\end{Con}
We will show that our main result, presented below, can be used to
infer that Beiter's conjecture
is `very false'.
\begin{Thm}
\label{main}
Let $p$ be a prime. Given an $1\le \beta\le p-1$ we let $\beta^*$ be the unique
integer $1\le \beta^*\le p-1$ with $\beta \beta^*\equiv 1({\rm mod~}p)$. \\
\indent Let
${\cal B}_{-}(p)$ be the set of integers $\beta$ satisfying
\begin{equation}
\label{verg2}
1\le \beta\le {p-3\over 2},~p\le \beta+2\beta^*+1,~\beta>\beta^*.
\end{equation}
For every prime $q\equiv \beta({\rm mod~}p)$ with 
$q>q_{-}(p)$ and
$\beta\in {\cal B}_{-}(p)$, there exists a prime $r_{-}>q$ and 
an integer $n_{-}$ such that $a_{pqr_{-}}(n_{-})=\beta_{-}-p$, where $q_{-}(p),r_{-}$ and $n_{-}$
can be explicitly given.\\
\indent Let
${\cal B}_+(p)$ be the set of integers $\beta$ satisfying
\begin{equation}
\label{verg1}
1\le \beta\le {p-3\over 2},~\beta+\beta^*\ge p,~\beta^*\le 2\beta,
\end{equation}
For every prime
$q\equiv \beta({\rm mod~}p)$ with $q>q_{+}(p)$ and $\beta \in 
{\cal B}_{+}(p)$ there exists a prime $r_{+}>q$ and 
an integer $n_{+}$ such that $a_{pqr_{+}}(n_{+})=p-\beta$, where $q_{+}(p),r_{+}$ and $n_{+}$
can be explicitly given. In case $\beta\in {\cal B}_{+}(p)$ and $\beta+\beta^*=p$,
then $A(pqr_{+})=p-\beta$.\\
\indent Put ${\cal B}(p)={\cal B}_{-}(p)\cup {\cal B}_{+}(p)$. If
${\cal B}(p)$ is non-empty, then
$$M(p)\ge p-{\min}\{{\cal B}(p)\}>{p+1\over 2},$$
and so Beiter's conjecture is false for the prime $p$.
\end{Thm}
Explicit choices of $q_{-}(p),r_{-}$ and $n_{-}$ are given in Theorem \ref{negative}
and explicit choices of $q_{+}(p),r_{+}$ and $n_{+}$ in Theorem \ref{positive}.\\
\indent Note that the sets ${\cal B}_{-}(p)$ and ${\cal B}_{+}(p)$ are disjoint.
For $p<11$ the set ${\cal B}(p)$ turns out to be empty.
For $11\le p\le 73$ it is given in Table 1. 
The underlined element is $(p-3)/2$ and for this range always turns out to
be in ${\cal B}(p)$.
The final column gives a lower
bound for $M(p)$. 
The table shows that Beiter's conjecture is false
for $11\le p\le 73$.
\begin{Prop}
\label{uno}
For $p\ge 11$, ${\cal B}(p)$ is non-empty and
$\max\{{\cal B}(p)\}={(p-3)/2}$.
\end{Prop}
{\it Proof}. Consider $\beta=(p-3)/2$. If $p\equiv 1({\rm mod~}3)$, then
$\beta^*=2(p-1)/3$ and one checks that $\beta\in {\cal B}_{+}(p)$. 
If $p\equiv 2({\rm mod~}3)$, then
$\beta^*=(p-2)/3$ and one checks that $\beta\in {\cal B}_{-}(p)$. \qed\\

\centerline{{\bf Table 1:} {\tt The sets ${\cal B}_{-}(p)$, ${\cal B}_{+}(p)$} and ${\cal B}(p)$}
\begin{center}
\begin{tabular}{|c|c|c|c|c|}
\hline
$p$ & ${\cal B}_{-}(p)$ & ${\cal B}_{+}(p)$ & ${\cal B}_{-}(p)\cup {\cal B}_{+}(p)={\cal B}(p)$ &
$p-\min{\cal B}(p)$ \\
\hline
11 & $\{{\underline 4}\}$ & $\emptyset$ & $\{4\}$ & 7\\
\hline
13 & $\emptyset$ & $\{{\underline 5} \}$ & $\{5\}$ & 8\\
\hline
17 & $\{{\underline 7}\}$ & $\emptyset $ & $\{7\}$ & 10\\
\hline
19 & $\emptyset$ & $\{{\underline 8} \}$ & $\{8\}$ & 11\\
\hline
23 & $\{{\underline{10}}\}$ & $\{9 \}$ & $\{9,10\}$ & 14\\
\hline
29 & $\{{\underline{13}}\}$ & $\{12 \}$ & $\{12,13\}$ & 17\\
\hline
31 & $\{13\}$ & $\{{\underline{14}} \}$ & $\{13,14\}$ & 18\\
\hline
37 & $\emptyset$ & $\{{\underline{17}} \}$ & $\{17\}$ & 20\\
\hline
41 & $\{18,{\underline{19}}\}$ & $\{17 \}$ & $\{17,18,19\}$ & 24\\
\hline
43 & $\{18\}$ & $\{19,\underline{20} \}$ & $\{18,19,20\}$ & 25\\
\hline
47 & $\{{\underline{22}}\}$ & $\{18,20 \}$ & $\{18,20,22\}$ & 29\\
\hline
53 & $\{{\underline{25}}\}$ & $\{22,23,24 \}$ & $\{22,23,24,25\}$ & 31\\
\hline
59 & $\{23,26,{\underline{28}}\}$ & $\{27 \}$ & $\{23,26,27,28\}$ & 36\\
\hline
61 & $\{25,28\}$ & $\{27,{\underline{29}} \}$ & $\{25,27,28,29\}$ & 36\\
\hline
67 & $\emptyset$ & $\{26,30,\underline{32} \}$ & $\{26,30,32\}$ & 40\\
\hline
71 & $\{32,33,{\underline{34}}\}$ & $\{27,29,30,31 \}$ & $\{27,29,30,31,32,33,34\}$ & 44\\
\hline
73 & $\{33\}$ & $\{27,30,34,{\underline{35}} \}$ & $\{27,30,33,34,35\}$ & 46\\
\hline
\end{tabular}
\end{center}

{}From Theorem \ref{main}, Proposition \ref{uno} and Dirichlet's theorem on
arithmetic progressions the following result is inferred.
\begin{Thm}
\label{main2} 
Suppose $p\ge 11$ is a prime.
Then the set of $q$ for which $A(pqr)>(p+1)/2$ for some prime $r$ has a positive
lower density ${\underline \delta}$ satisfying
$${\underline \delta}\ge {|{\cal B}_{+}(p)|+|{\cal B}_{-}(p)|\over p-1}\ge {1\over p-1}>0,$$
and hence Beiter's conjecture is false for the prime $p$.
\end{Thm}
The elementary method of proof of Proposition 1 allows one to prove, e.g., that if
$p\equiv 23({\rm mod~}24)$ and $p>23$, then $M(p)\ge (5p-3)/8$ and more generally
it allows one to indicate for every $\epsilon>0$ an arithmetic progression such
that for all primes larger than some explicit number in this progression we
have $M(p)\ge ({2\over 3}-\epsilon)p$ (Proposition \ref{yves}). On invoking a result from 
the theory of inverses modulo $p$ it can be even shown that the latter lower bound holds
for {\it all} primes $p$ sufficiently large.
\begin{Thm}
\label{main3} 
Let $\epsilon>0$. Then ${2\over 3}p(1-\epsilon)\le M(p)\le {3\over 4}p$ for every sufficiently 
large prime $p$.
\end{Thm}
In Table 2 we give for some small primes $p$ intervals $[a,b]$ such that
$a\le M(p)\le b$. The number $b=p-\lceil p/4\rceil$ and $a=|a_{pqr}(n)|$, showing
that $M(p)\ge a$. The values of $a$ are the largest known to us and were found
by extensive computer calculation.\\

\centerline{{\bf Table 2:} {\tt Interval for $M(p)$}}
\begin{center}
\begin{tabular}{|c|c|c|c|c|c|c|c|}
\hline
$p$ & $(p+1)/2$ & $q$ & $r$ & $n$ & $M(p)$ {\rm interval} & $[2p/3]$ &
$\delta(p)$\\
\hline
3 & 2 & 5 & 7 & 7 & [2, 2] & 2 & 0\\
\hline
5 & 3 & 7 & 11 & 119 & [3, 3] & 3 & 0\\
\hline
7 & 4 & 11 & 37 & 963 & [4, 5] & 4 & $\ge 0$\\
\hline
11 & 6 & 19 & 601 & 34884 & [7, 8] & 7 & $\ge 1$\\
\hline
13 & 7 & 31 & 1097 & 137160 & [8, 9] & 8 & $\ge 1$\\
\hline
17 & 9 & 29 & 41 & 4801 & [10, 12] & 11 & $\ge 1$\\
\hline
19 & 10 & 53 & 859 & 318742 & [12, 14] & 12 & $\ge 2$\\
\hline
23 & 12 & 41 & 4903 & 1583731 & [14, 17] & 15 & $\ge 2$\\
\hline
29 & 15 & 127 & 7793 & 8915220 & [18, 21] & 19 & $\ge 3$\\
\hline
31 & 16 & 89 & 4519 & 4424131 & [19, 23] & 20 & $\ge 3$\\
\hline
37 & 19 & 47 & 1217 & 743670 & [22, 27] & 24 & $\ge 3$\\
\hline
41 & 21 & 71 & 97 & 96529 & [26, 30] & 27 & $\ge 5$\\
\hline
43 & 22 & 53 & 2963 & 2358548 & [26, 32] & 28 & $\ge 4$\\
\hline
47 & 24 & 347 & 12113 & 64756445 & [29, 35] & 31 & $\ge 5$\\
\hline
53 & 27 & 61 & 17377 & 18037438 & [33, 39] & 35 & $\ge 6$\\
\hline
59 & 30 & 67 & 21247 & 27047555 & [37, 44] & 39 & $\ge 7$\\
\hline
61 & 31 & 191 & 30203 & 126913006 & [38, 45] & 40 & $\ge 7$\\
\hline
67 & 34 & 191 & 91127 & 417817361 & [42, 50] & 44 & $\ge 8$\\
\hline
71 & 36 & 311 & 13327 & 91183645 & [44, 53] & 47 & $\ge 8$\\
\hline
73 & 37 & 83 & 4241 & 9156474 & [46, 54] & 48 & $\ge 9$\\
\hline
\end{tabular}
\end{center}
The last column gives information
about the difference $\delta(p):=M(p)-(p+1)/2$. In case $\delta(p)>0$
the associated $p,q$ and $r$ give rise to a counter-example to
Beiter's conjecture.\\ 
\indent If $\beta\in {\cal B}_{-}(p)$,
then $p\le 3\beta-1$ and $p-\beta\le (2p-1)/3$. 
If $\beta\in {\cal B}_{+}(p)$, then $p\le \beta+\beta^*\le 3\beta$
and so $\beta\ge p/3$ and hence $p-\beta\le 2p/3$.
Thus Theorem \ref{main} only
allows one to find counter-examples $\le 2p/3$ to Beiter's conjecture. Extensive
numerical computations gave many counter-examples not covered by Theorem \ref{main},
but all of them are $\le 2p/3$.
Thus the strongest corrected version
of Beiter's conjecture which we can presently neither disprove nor prove is as follows.
\begin{Con} {\rm (Corrected Beiter conjecture)}.
We have $M(p)\le 2p/3$.
\end{Con}
Note that it implies that Beiter's original conjecture is correct for $p=7$. This is
at present still an open problem. If $A(7qr)>4$, then we must have $q\equiv \pm 3({\rm mod~}7)$ by
Theorem \ref{flup}.\\
\indent Our final result deals with some apparent variations of $M(p)$.
Let $M_{+}(p)$ and $M_{-}(p)$ be the maximum, respectively
minimum of the coefficients of the ternary cyclotomic polynomials
with $p$ the smallest prime factor of $n$.
\begin{Thm}
\label{plus}
We have $M_{-}(p)=M_{+}(p)=M(p)$.
\end{Thm}

\indent For a nice survey of properties of coefficients of cyclotomic polynomials
see Thangadurai \cite{Thanga}.

\subsection{Some results of Nathan Kaplan}
Using the identity
$$\Phi_{pqr}(x)=(1+x^{pq}+x^{2pq}+\cdots)(1+x+\cdots+x^{p-1}-x^q-\cdots-x^{q+p-1})
\Phi_{pq}(x^r),$$
Kaplan \cite{Kaplan} proved the following lemma.
\begin{Lem} {\rm (Nathan Kaplan, 2007)}.
\label{kapel}
Let $n\ge 0$ be an integer.
Put $$b_i=\cases{a_{pq}(f(i)) & if $f(i)\le n/r$; \cr
0 & otherwise,}$$ 
where $f(m)$ is the unique value $0\le  f(m) < pq$ such that
$$f(m)\equiv {n-m\over r} ~({\rm mod~}pq).$$
Then
\begin{equation}
\label{lacheens}
a_{pqr}(n)=\sum_{m=0}^{p-1}b_{m}-\sum_{m=0}^{p-1}b_{m+q},
\end{equation}
\end{Lem}
\noindent Since the $a_{pq}(i)$ are easily computed, Kaplan's lemma is actually
useful. Indeed, his lemma plays a crucial role
in our counter-example constructions. 
A nice feature of the lemma is that it works for every $n\ge 0$. Thus if it
shows that $a_{pqr}(n)\ne 0$, then we know that $n\le \varphi(pqr)$. In our
counter-example constructions this saves us from checking that for the chosen
$n$ we have $n\le \varphi(pqr)$.\\
\indent The next lemma gives the values of $a_{pq}(i)$. For a proof
see e.g. Lam and Leung \cite{LL} or Thangadurai \cite{Thanga},
\begin{Lem}
\label{binary}
Let $p<q$ be odd primes. Let $\rho$ and $\sigma$ be the (unique) non-negative
integers for which $(p-1)(q-1)=\rho p+\sigma q$
Let $0\le m<pq$. Then either $m=\alpha_1p+\beta_1q$ or $m=\alpha_1p+\beta_1q-pq$
with $0\le \alpha_1\le q-1$ the unique integer such that $\alpha_1 p\equiv m({\rm mod~}q)$
and $0\le \beta_1\le p-1$ the unique integer such that $\beta_1 q\equiv m({\rm mod~}p)$.
The cyclotomic coefficient $a_{pq}(m)$ equals
$$\cases{1 & if $m=\alpha_1p+\beta_1q$ with $0\le \alpha_1\le \rho,~0\le \beta_1\le
\sigma$;\cr -1 & if $m=\alpha_1p+\beta_1q-pq$ with $\rho+1\le \alpha_1\le q-1,~\sigma+1\le 
\beta_1\le p-1$;\cr  0 & otherwise.}
$$
\end{Lem}
The following result of Kaplan \cite{Kaplan} together
with Dirichlet's theorem shows that given one counter-example
$(p,q,r)$ infinitely
many counter-examples to Beiter's conjecture exist with
the same values of $p$ and $q$.
\begin{Thm} {\rm (Nathan Kaplan, 2007)}. \label{kapel2}
For any prime $s>q$ such that $s\equiv \pm r({\rm mod~}pq)$ we have
$A(pqr)=A(pqs)$.
\end{Thm} 
{\tt Example}. Put $p=17$ and $q=29$. By computation one finds that
$A(pq\cdot 1931)=10$. On applying Kaplan's result one then finds
from this that also $A(pq\cdot 2917)=A(pq\cdot 2999)=10$.\\

\noindent Implicit in Kaplan's proof of Theorem \ref{kapel2} is the following
result using which we immediately infer that Theorem \ref{plus} holds true.
\begin{Prop}
Suppose that $a_{pqr}(n)=m$. Write
$n=[{n\over r}]r+n_0$ with $0\le n_0<r$.\\
{\rm 1)} Let $s>r$ be a prime satisfying $s\equiv r({\rm mod~}pq)$. Then
$$a_{pqs}\Big(\Big[{n\over r}\Big]s+n_0\Big)=m.$$
{\rm 2)} Let $t>pq$ be a prime satisfying $t\equiv -r({\rm mod~}pq)$. Let
$0\le n_1<pq$ be the unique integer such that 
$n_1\equiv q+p-1-n_0({\rm mod~}pq)$. Then
$$a_{pqt}\Big(\Big[{n\over r}\Big]t+n_1\Big)=-m.$$
\end{Prop}

\section{A counter-example construction for $p=11$}
Using Lemma \ref{kapel} and Lemma \ref{binary} (we leave this
as an exercise to the reader) one finds that for the $p,q$ and $r$ in the M\"oller
construction we have $b_{f(m)}=1$ for $0\le m\le (p-1)/2$ and 
$b_{f(m)}=0$ for the remaining $m$ in (\ref{lacheens}), giving
$a_{pqr}(n)=(p+1)/2$. Likewise,
for the Lehmer example we find $b_{f(m)}=1$ for $0\le m\le (p-3)/2$ and 
$b_{f(m)}=0$ for the remaining $m$, giving $a_{pqr}(n)=(p-1)/2$.\\ 
\indent For the counter-examples to Beiter's conjecture we find
in general rather more complicated vectors $(b_{f(0)},\ldots,b_{f(p-1)})$ and 
$(b_{f(q)},\ldots,b_{f(q+p-1)})$. However, some of them are regular enough
as to build a general construction on. We give
an example which is intended as an appetizer that should help the reader digest more
easily the general construction given in Theorem \ref{negative}. Notice that
Theorem \ref{negative} implies Theorem \ref{11}. The first few counter-examples
produced by Theorem \ref{11} are given in Table 3.\\
\vfil\eject
\centerline{{\bf Table 3:} {\tt Some counter-examples produced by Theorem \ref{11}}}
\begin{center}
\begin{tabular}{|c|c|c|c|c|c|}
\hline
$p$ & $q$    & $\alpha$ & $r$ & $n$  & $a_{pqr}(n)$ \\
\hline
11 & 59 & 2 & 877 & 175410 & -7\\
\hline
 & 103 & 4 & 1229 & 381000 & -7\\
\hline
 & 191 & 6 & 4639 & 3173086 & -7\\
\hline
 &  & 7 & 16937 & 10280769 & -7\\
\hline
 & 257 & 8 & 3011 & 2788196 & -7\\
\hline
 &  & 9 & 1163 & 987397 & -7\\
\hline
 & & 10 & 8731 & 6740342 & -7\\
\hline
 & 367 & 12 & 56999 & 72844732 & -7\\
\hline
 & & 13 & 811 & 974021 & -7\\
\hline
 & & 14 & 39157 & 44012478 & -7\\
\hline
\end{tabular}
\end{center}

\begin{Thm}
\label{11}
Let $q<r$ be primes such that $q\equiv 4({\rm mod~}11)$
and $r\equiv   -3({\rm mod~}11)$. Let
$1\le \alpha\le q-1$ be the unique integer such that
$11r\alpha\equiv 1({\rm mod~}q)$.
Suppose that
$${q\over 33}< \alpha\le {3q-1\over 77}.$$
Then $a_{11qr}(10+(6q-77\alpha)r)=-7$.
\end{Thm}
{\it Proof}. Put $p=11$ and $n=10+(6q-77\alpha)r$. 
Note that $n\ge 10$. We will compute $a_{pqr}(n)$ using Lemma \ref{kapel}. Since this
will have -7 as outcome, it follows that $n\le \varphi(pqr)$.
Observe that 
$$0\le {n-10\over r}=(q-7\alpha)p+6q-pq\le 6q<pq$$ and so $f(10)=(q-7\alpha)p+6q-pq$. 
This expresses $f(10)$ as a linear combination in $p$ and $q$. 
Let $0\le r_1<pq$ be the unique integer with
$r_1\equiv -{1\over r}({\rm mod~}pq)$. It is easy to see
that $r_1=q-\alpha p$.
On
noting that $f(m)\equiv f(10)+(m-10)r_1\equiv f(10)+(m-10)(q-\alpha p)~({\rm mod~}pq)$,
we infer that $f(m)$
is congruent modulo $pq$ to the corresponding entry in Table 4 (and likewise
for $f(m+q)$ on 
noting that 
$-q/r\equiv 4q({\rm mod~}pq)$ and $f(m+q)\equiv f(m)-q/r\equiv f(m)+4q({\rm mod~}pq)$). 
In the $f(m+q)$ column we trivially have
$$0\le f(q)<\cdots<f(q+p-1)=10q-7\alpha p\le 10q<pq.$$
Using this we infer that $f(m+q)$ is actually equal to the
corresponding entry in Table 4. In the $f(m)$ column we 
have $0\le f(0)<\cdots<f(4)<pq$ and on using that $\alpha\le q/(2p)$ we
find $0\le f(5)<\cdots<f(10)<6q<pq$. Again we see that $f(m)$ is 
actually equal to the corresponding entry in Table 4.\\
\indent Now we are ready to invoke Lemma \ref{binary}. One computes that
$\sigma=2$ and $\rho=(8q-10)/11$. The conditions on $\alpha$ ensure that
$q-7\alpha\ge \rho+1$ and $3\alpha\le \rho$. On applying Lemma \ref{binary} we
then infer that $c_m:=a_{pq}(f(m))$ and $c_{m+q}:=a_{pq}(f(m+q))$ are as 
given in Table 4. Since $[n/r]=f(10)$ it follows, by Lemma \ref{kapel} that
$b_m=c_m$ if $f(m)\le f(10)$ and $b_m=0$ otherwise. Thus to compute say
the $b_m$ column we have $b_m=0$ if $c_m=0$. If $c_m\ne 0$ we have
$$b_m=\cases{c_m & if $f(m)\le f(10)$;\cr
0 & otherwise.}$$
Note that clearly $f(7)<....<f(10)$. It then follows that the $b_m$ column
equals the $c_m$ column. Next let us determine the $c_{m+q}$ column.
We claim that $f(q)<f(q+1)<f(q+2)<f(10)$. To establish this we have to check that 
$\alpha p+2q<(q-7\alpha)p-5q$. Note that $f(q+4)<\cdots<f(q+10)$.
The conditions on $\alpha$ ensure that $f(q+4)>f(10)$ and we
see that $f(10)<f(q+4)<\cdots<f(q+10)$ and thus $b_{q+4}=\ldots=b_{q+10}=0$.
Finally, on applying Lemma \ref{kapel} we infer that
$$a_{pqr}(n)=\sum_{m=0}^{p-1}b_{m}-\sum_{m=0}^{p-1}b_{m+q}=-4-3=-7.$$
This completes the proof. \qed\\

\noindent {\tt Remark}. Note that the conditions imposed on $\alpha$ are such
that $\alpha\le q/(2p)$, $3\alpha\le \rho$, $q-7\alpha\ge \rho+1$ (which is 
equivalent with $\alpha\le (3q-1)/77$),
$f(q+4)>f(10)$ (which is equivalent with $\alpha>q/33$) and $f(q+2)\le f(10)$.\\

\centerline{{\bf Table 4:} {\tt A counter-example construction for $p=11$}}
\begin{center}
\begin{tabular}{|c|c|c|c|c|c|c|c|}
\hline
$m$ & $f(m)$    & $c_m$ & $b_m$ & $m+q$ & $f(m+q)$ & $c_{m+q}$ & $b_{m+q}$\\
\hline  
0  & $3\alpha p+q(p-4)$ & 0 & 0 & $q$ & $3\alpha p$ & 1 & 1\\
\hline  
1 & $2\alpha p +q(p-3)$ & 0 & 0 & $q+1$ &  $2\alpha p+q$ & 1 & 1\\
\hline
2 & $\alpha p +q(p-2)$ & 0 & 0 & $q+2$ &  $\alpha p+2q$ & 1 & 1\\
\hline
3 & $q(p-1)$ & 0 & 0 & $q+3$ &  $3q$ & 0 & 0\\
\hline
4 & $(q-\alpha)p$ & 0 & 0 & $q+4$ &  $(q-\alpha)p+4q-pq$ & -1 & 0\\
\hline
5 & $(q-2\alpha)p +q-pq$ & 0 & 0 & $q+5$ &  $(q-2\alpha)p+5q-pq$ & -1 & 0\\
\hline
6 & $(q-3\alpha)p +2q-pq$ & 0 & 0 & $q+6$ &  $(q-3\alpha)p+6q-pq$ & -1 & 0\\
\hline
7 & $(q-4\alpha)p +3q-pq$ & -1 & -1 & $q+7$ &  $(q-4\alpha)p+7q-pq$ & -1& 0\\
\hline
8 & $(q-5\alpha)p +4q-pq$ & -1 & -1 & $q+8$ &  $(q-5\alpha)p+8q-pq$ & -1 & 0\\
\hline
9 & $(q-6\alpha)p +5q-pq$ & -1 & -1 & $q+9$ &  $(q-6\alpha)p+9q-pq$ & -1 & 0\\
\hline
10 & $(q-7\alpha)p +6q-pq$ & -1 & -1 & $q+10$ &  $(q-7\alpha)p+10q-pq$ & -1 & 0\\
\hline
\end{tabular}
\end{center}

\section{General counter-example construction}
\subsection{The negative coefficient case}
We now establish a more general counter-example construction. The approach
will be similar to that of the previous section. For reasons of space
the analogue of Table 4, Table 5, is split into two tables, for $f(m)$,
respectively $f(m+q)$.\\
\vfil\eject
\centerline{{\bf Table 5A:} {\tt General negative coefficient construction, $f(m)$ case}}
\begin{center}
\begin{tabular}{|c|c|c|c|}
\hline
$m$ & $f(m)$ & $c_{m}$ & $b_{m}$\\
\hline
0 & $(\sigma+1)\alpha p+(p-\beta)q$ & 0 & 0\\
\hline
1 & $\sigma \alpha p+(p-\beta+1)q$ & 0 & 0\\
\hline
$\cdots$ & $\cdots$ & 0 & 0\\
\hline
$\sigma+1$ & $0\cdot p+(\sigma+1+p-\beta)q$ & 0 & 0\\
\hline
$\sigma+2$ & $(q-\alpha)p+(\sigma+2+p-\beta)q-pq$ & -1 & 0\\
\hline
$\cdots$ & $\cdots$ & -1 & 0\\
\hline
$\beta-1$ & $(q-(\beta-\sigma-2)\alpha)p+(p-1)q-pq$ & -1 & 0\\
\hline
$\beta$ & $(q-(\beta-\sigma-1)\alpha)p+0\cdot q$ & 0 & 0\\
\hline
$\beta+1$ & $(q-(\beta-\sigma)\alpha)p+1\cdot q-\delta_1 pq$ & 0 & 0\\
\hline
$\cdots$ & $\cdots$ & 0 & 0\\
\hline
$\beta+\sigma$ & $(q-(\beta-1)\alpha)p+\sigma q-\delta_{\sigma}pq$ & 0 & 0\\
\hline
$\beta+\sigma+1$ & $(q-\beta\alpha)p+(\sigma+1) q-pq$ & -1 & -1\\
\hline
$\cdots$ & $\cdots$ & -1 & -1\\
\hline
$\beta+k$ & $(q-(\beta-\sigma-1+k)\alpha)p+kq-pq$ & -1 & -1\\
\hline
$\cdots$ & $\cdots$ & -1 & -1\\
\hline
$p-1$ & $(q-(p-\sigma-2)\alpha)p+(p-\beta-1)q-pq$ & -1 & -1\\
\hline
\end{tabular}
\end{center}
In Table 5A $\delta_j$ is the unique integer such that the corresponding
entry is in the interval $[0,pq)$.\\

\centerline{{\bf Table 5B:} {\tt General negative coefficient construction, $f(m+q)$ case}}
\begin{center}
\begin{tabular}{|c|c|c|c|}
\hline
$m+q$ & $f(m+q)$ & $c_{m+q}$ & $b_{m+q}$\\
\hline
$q$ & $(\sigma+1)\alpha p$ & 1 & 1\\
\hline
$q+1$ & $\sigma \alpha p+q$ & 1 & 1\\
\hline
$\cdots$ & $\cdots$ & 1 & 1\\
\hline
$q+\sigma$ & $\alpha p+\sigma q$ & 1 & 1\\
\hline
$q+\sigma+1$ & $0\cdot p+(\sigma+1)q$ & 0 & 0\\
\hline
$q+\sigma+2$ & $(q-\alpha)p+(\sigma+2)q-pq$ & -1 & 0\\
\hline
$\cdots$ & $\cdots$ & -1 & 0\\
\hline
$q+p-1$ & $(q-(p-\sigma-2)\alpha)p+(p-1)q-pq$ & -1 & 0\\
\hline
\end{tabular}
\end{center}

\begin{Lem}
\label{mastercopy}
Let $p$ be a prime. Let $1\le \beta\le (p-3)/2$. Let $q>p$ be a prime
satisfying $q\equiv \beta({\rm mod~}p)$ and $r>q$ be a prime
satisfying $qr\equiv -1({\rm mod~}p)$. Let $1\le \alpha\le q-1$
be the unique integer such that $pr\alpha\equiv 1({\rm mod~}q)$.
Put $$w_{-}=(p-\beta-1)q-(p-\sigma-2)\alpha p,$$ where $\rho$ and $\sigma$
are uniquely determined by $(p-1)(q-1)=\rho p+\sigma q,~\rho,\sigma\ge 0$.
Suppose that
$$p\ge \beta+\sigma+2,~\beta\ge \sigma+2$$
and
\begin{equation}
\label{a}
\alpha\le {q(\sigma+1)\over p\beta},
\end{equation}
\begin{equation}
\label{b}
\alpha\le {q(p-1-\sigma)-(p-1)\over p(\sigma+1)},
\end{equation}
\begin{equation}
\label{c}
\alpha\le {q(\sigma+1)-1\over p(p-\sigma-2)},
\end{equation}
\begin{equation}
\label{d}
\alpha\le {q(p-\sigma-1-\beta)\over p(p-\sigma-1)},
\end{equation}
and
\begin{equation}
\label{e}
\alpha> {q(p-3-\sigma-\beta)\over p(p-\sigma-3)},
\end{equation}
then $a_{pqr}(p-1+rw_{-})=\beta-p$.
\end{Lem}
{\tt Remark}. The conditions (\ref{a}), (\ref{b}), (\ref{c}), (\ref{d}) and
(\ref{e}) are used to ensure that respectively,
$f(\beta+\sigma+1)\ge 0$, $(\sigma+1)\alpha\le \rho$, $q-(p-\sigma-2)\alpha\ge \rho+1$,
$f(q+\sigma)\le f(p-1)$
and $f(q+\sigma+2)>f(p-1)$.\\ 

\noindent {\it Proof of Lemma} \ref{mastercopy}. 
Let $0\le r_1<pq$ be the unique integer with
$r_1\equiv -{1\over r}({\rm mod~}pq)$. 
We have $r_1\equiv q-\alpha p ({\rm mod~}pq)$.
By (\ref{a}) and since $\beta\ge \sigma+2$, we infer that $q-\alpha p=r_1$.
Reasoning as in the proof of Theorem \ref{11} we find
that $f(m)$ and $f(m+q)$ are congruent modulo $pq$ to the numbers given in
Table 5A, respectively 5B.\\
\indent In Table 5A we distinguish 3 ranges: $0\le m\le \beta$, $\beta+1\le m\le \beta+
\sigma$ and $\beta+\sigma+1\le m\le p-1$. In the first range the entries 0 up to $\beta$
are non-negative and in ascending order. Since the entry for $\beta$ is $<pq$
it follows that for $0\le m\le \beta$, $f(m)$ is actually equal to the corresponding
entry given in Table 5A. Since $(\sigma+1)\alpha\le \rho$ and $p-\beta\ge \sigma+1$
it follows by Lemma \ref{binary} that $c_m=0$ (and hence $b_m=0$) for $m=0,\ldots,\sigma+1$.
For $\sigma+2\le m\le \beta-1$ one finds that $c_m=-1$. Since
$$f(\sigma+2)=f(q+\sigma+2)+(p-\beta)q>f(q+\sigma+2)>f(p-1),$$
we infer that $b_m=0$ for $\sigma+2\le m\le \beta-1$. Clearly $b_{\beta}=c_{\beta}=0$.\\
\indent In the second range we have, for $1\le k\le \sigma$, 
$$(q-(\beta-\sigma+k-1)\alpha)p+kq-\delta_kpq$$
as entry in row $\beta+k$, where a priori $\delta_k$ is an integer.
Since
$$q-1\ge q-(\beta-\sigma+k-1)\alpha\ge q-(p-\sigma-2)\alpha\ge \rho+1>0,$$
we find that $0\le (q-(\beta-\sigma+k-1)\alpha)p+kq<2pq$ and hence
$\delta_k\in \{0,1\}$. By Lemma \ref{binary} again we now find that
$c_m=b_m=0$ for $\beta+1\le m\le \beta+\sigma$.\\
\indent In the range $\beta+\sigma+1\le m\le p-1$
the entries in Table 5A are in ascending order and the final entry is less
than $pq$. Since $\alpha\le (\sigma+1)q/(\beta p)$, it follows that the
$\beta+\sigma+1$ entry in the table is $\ge 0$. It follows that $f(m)$ is actually
equal to the corresponding entry in Table 5A. 
Since $q-(p-\sigma-2)\alpha \ge \rho+1$
 it now follows By Lemma \ref{binary} that
$c_m=-1$ for $\beta+\sigma+1\le m\le p-1$. Since $f(m)\le f(p-1)$ for
$\beta+\sigma+1\le m\le p-1$ we infer that also $b_m=-1$ in this range.\\
\indent Establishing the correctness of the $f(m+q)$ and $c_{m+q}$ column is
straightforward and left to the reader. Note that once we have
$f(q+\sigma)\le f(p-1)$ and $f(q+\sigma+2)>f(p-1)$, the $b_{m+q}$ column is
as given in the table. That these inequalities hold is ensured by conditions
(\ref{d}), respectively (\ref{e}). 
Finally on applying Lemma \ref{kapel} we infer that
$$a_{pqr}(p-1+rw_{-})=-\sum_{j=\beta+\sigma+1}^{p-1}1-\sum_{j=q}^{q+\sigma}1
=-(p-1-\beta-\sigma)-(\sigma+1)=\beta-p.$$
This concludes the proof.\qed\\

\noindent In the next lemma the set of {\it real} numbers $\alpha$ satisfying
(\ref{a}), (\ref{b}), (\ref{c}), (\ref{d}) and (\ref{e}) is determined.
\begin{Lem}
\label{abcde}
Let ${\cal I}$ be the set of real numbers satisfying {\rm (\ref{a}), 
(\ref{b}), (\ref{c}), (\ref{d}) and (\ref{e})} and suppose the conditions
of Lemma {\rm \ref{mastercopy}} preceding {\rm (\ref{a})} are satisfied. 
Then the set ${\cal I}$
is non-empty iff $p\le \beta+2\beta^*+1$. In that case
$${\cal I}=\cases{\Big({q\over p}\Big(1-{\beta\over p-\beta^*-2}\Big),
{q\over p}\Big(1-{\beta\over p-\beta^*}\Big)\Big] & if $p<\beta+2\beta^*+1$;\cr
\Big({q\over p}\Big(1-{\beta\over p-\beta^*-2}\Big),
{q\beta^*-1\over p(p-\beta^*-1)}\Big] & if $p=\beta+2\beta^*+1$.}$$
If ${\cal I}$ is non-empty then it consists of positive reals only.
\end{Lem}
{\it Proof}. {}From $(\rho+1)p+(\sigma+1)q=pq+1$ we infer that
$(\sigma+1)q\equiv 1({\rm mod~}p)$ and hence $\sigma\equiv 1/q-1({\rm mod~}p)$.
Since $q\equiv \beta({\rm mod~}p)$, we infer that $\sigma=\beta^*-1$. Note
that if (\ref{c}) is satisfied then, since $p-\sigma-2\ge \beta$, automatically
(\ref{a}) is satisfied.
Note that if $\alpha$ satisfies (\ref{c}),
then 
$$\alpha\le {q(\sigma+1)\over p(p-\sigma-3)}.$$
Now if $\alpha$ is also to satisfy (\ref{e}), then we must have
$\sigma+1>p-3-\sigma-\beta$ and so $p\le 2\sigma+3+\beta$. 
Thus if $p>\beta+2\beta^*+1$, then ${\cal I}$ is empty and hence we may assume
that $p\le \beta+2\beta^*+1$.  
By (\ref{d}) and since $p\ge 2\sigma+3=2\beta^*+1$ we infer that
$$\alpha\le {q\over p}\le {q(p-2-\sigma)\over p(\sigma+1)}\le {q(p-1-\sigma)-(p-1)\over
p(\sigma+1)},$$
and hence the condition (\ref{b}) is also superfluous. Note that the
inequality
$${q(p-\sigma-1-\beta)\over p(p-\sigma-1)}\le {q(\sigma+1)-1\over p(p-\sigma-2)}$$
can be rewritten as
$$p\le \beta+3+2\sigma-{1\over q}-{\beta\over p-\sigma-1}.$$
We observe that
$$0<{1\over q}+{\beta\over p-\sigma-1}\le {1\over q}+{\beta\over \beta+1}\le {1\over q}+1-{1\over
p}<1.$$
It follows that if $p=\beta+3+2\sigma=\beta+2\beta^*+1$, then (\ref{d})
is redundant and if $p<\beta+2\beta^*+1$, then (\ref{c}) is redundant. In the
latter case we obtain that 
$${\cal I}=\Big({q\over p}\Big(1-{\beta\over p-\beta^*-2}\Big),
{q\over p}\Big(1-{\beta\over p-\beta^*}\Big)\Big],$$
a clearly non-empty interval. In the former case we obtain
$${\cal I}=\Big({q\over p}\Big(1-{\beta\over p-\beta^*-2}\Big),
{q\beta^*-1\over p(p-\beta^{*}-1)}\Big)\Big],$$
in which case an easy calculation shows that it is a non-empty interval.
Since $\beta+\beta^*\le 2\beta-1\le p-4$ it follows that $p-\beta^*-2>\beta$ and
thus if ${\cal I}$ is non-empty, it contains positive reals only. \qed

\begin{Thm}
\label{negative}
Suppose that ${\cal B}_{-}(p)$ is non-empty and $\beta \in {\cal B}_{-}(p)$.
Suppose also that $p<\beta+2\beta^*+1$.
Let $q>p$ be a prime satisfying $q\equiv \beta({\rm mod~}p)$ and
$q\ge q_{-}(p)$ with
$q_{-}(p)=p(p-\beta^*)(p-\beta^*-2)/(2\beta)$. Then the interval
$${\cal I}=\Big({q\over p}\Big(1-{\beta\over p-\beta^*-2}\Big),
{q\over p}\Big(1-{\beta\over p-\beta^*}\Big)\Big]$$
contains at least one integer $a$.
Let $r>q$ be a prime with
$r(q-pa)\equiv -1({\rm mod~}pq)$,
then
$$a_{pqr}(p-1+[(p-\beta-1)q-(p-\beta^*-1)ap]r)=\beta-p<-{(p+1)\over 2}$$
is a counter-example to Beiter's conjecture.\\
\indent In case $p=\beta+2\beta^*+1$ the same conclusion holds, 
but with $q_{-}(p)$ replaced by $(\beta+\beta^*-1)(p(\beta+\beta^*)+1)/\beta$
and ${\cal I}$ by 
$${\cal I}=\Big({q\over p}\Big(1-{\beta\over p-\beta^*-2}\Big),
{q\beta^*-1\over p(p-\beta^*-1)}\Big].$$
\end{Thm}
{\it Proof}. Note that as a function of $q$ the length of ${\cal I}$ is
increasing. If the length of ${\cal I}$ is at least one, then it contains
at least one positive integer. Now $q_{-}(p)$ is obtained on solving
the equation $|{\cal I}|=1$ for $q$. The proof is completed on invoking 
Lemma \ref{mastercopy} and Lemma \ref{abcde}. \qed\\

\noindent Remark. If $1\le \beta\le (p-3)/2$ and $\beta>\beta^*$, then
$\beta$ satisfies the conditions of Lemma \ref{mastercopy}. If, in addition, $p\le \beta+2\beta^*+1$, that
is if $\beta\in {\cal B}_{-}(p)$, then the conditions of 
both Lemma \ref{mastercopy} and Lemma \ref{abcde} are satisfied by $\beta$.

\subsection{The positive coefficient construction}
Since the method of proof in this section is similar to that in the previous section, 
some of the details will be suppressed.
\begin{Lem}
\label{main1}
Let $p$ be a prime. Let $1\le \beta\le (p-3)/2$. Let $q>p$ be a prime
satisfying $q\equiv \beta({\rm mod~}p)$ and $r>q$ be a prime
satisfying $qr\equiv -1({\rm mod~}p)$. Let $1\le \alpha\le q-1$
be the unique integer such that $pr\alpha\equiv 1({\rm mod~}q)$.
Put $$w_{+}=1+(p-\beta-1)q-(p-\beta)\alpha p,$$ where $\rho$ and $\sigma$
are uniquely determined by $(p-1)(q-1)=\rho p+\sigma q,~\rho,\sigma\ge 0$.
Suppose that
$$\beta+\sigma\ge p-1$$
and
\begin{equation}
\label{a1}
\alpha< {(\sigma+1)q-1\over p(\beta-1)},
\end{equation}
\begin{equation}
\label{b1}
\alpha\le {q(p-\sigma-1)\over p(p-\beta)},
\end{equation}
\begin{equation}
\label{c1}
\alpha\le {q(\sigma+1-\beta)\over p(\sigma+1)},
\end{equation}
and
\begin{equation}
\label{d1}
\alpha> {q(p-2\beta-1)\over p(p-\beta-1)},
\end{equation}
then $a_{pqr}(p-1+rw_{+})=p-\beta$.
\end{Lem}

\noindent {\tt Remark}. The conditions (\ref{a1}), (\ref{b1}), (\ref{c1}) and (\ref{d1}) 
are used to ensure that respectively,
$\rho+1+(\beta-1)\alpha\le q-1$, $\rho+1-(p-\beta)\alpha\ge 0$, $f(q+p-\sigma-2)\le f(p-1)$
and $f(q+\beta)>f(p-1)$.\\

\noindent {\it Proof of Lemma} \ref{main1}. The proof is a more
general variant of the proof of Theorem \ref{11}. Again we make
use of tables in the proof. \\

\centerline{{\bf Table 6A:} {\tt General positive coefficient construction: $f(m)$ case}}
\begin{center}
\begin{tabular}{|c|c|c|c|c|c|}
\hline
$m$ & $f(m)$   & $c_m$ & $b_m$ \\
\hline  
0  & $(\rho+1+(\beta-1)\alpha)p + (\sigma-\beta+1)q-\tau_0pq$ & 0 & 0 \\
\hline  
$\cdots$ & $\cdots$ & 0 & 0 \\
\hline
$\beta-2$ & $(\rho+1+\alpha ) p +(\sigma-1)q-\tau_{\beta-2}pq$ & 0 & 0  \\
\hline
$\beta-1$ & $(\rho+1) p +\sigma q-\tau_{\beta-1}pq$ & 0 & 0  \\
\hline
$\beta$ & $(\rho+1-\alpha) p +(\sigma +1) q-\tau_{\beta}pq$ & 0 & 0  \\
\hline
$\cdots$ & $\cdots$ & 0 & 0 \\
\hline
$p-\sigma+\beta-2$ & $(\rho+1-(p-\sigma-1)\alpha)p+(p-1)q-\tau_{p-\sigma+\beta-2}pq$ & 0 & 0 \\
\hline
$p-\sigma+\beta-1$ & $(\rho+1-(p-\sigma)\alpha)p$ & 1 & 1 \\
\hline
$p-\sigma+\beta$ & $(\rho+1-(p-\sigma+1)\alpha)p+q$ & 1 & 1 \\
\hline
$\cdots$ & $\cdots$ & 1 & 1 \\
\hline
$p-1$ & $(\rho+1-(p-\beta)\alpha)p+(\sigma-\beta)q$ & 1 & 1 \\
\hline
\end{tabular}
\end{center}
In Table 6A $\tau_j$ is the unique integer in $\{0,1\}$ such that the corresponding
entry for $f(m)$ is in the interval $[0,pq)$.\\
\vfil\eject

\centerline{{\bf Table 6B:} {\tt General positive coefficient construction: $f(m+q)$ case}}
\begin{center}
\begin{tabular}{|c|c|c|c|c|}
\hline
$m+q$ & $f(m+q)$    & $c_{m+q}$ & $b_{m+q}$ \\
\hline  
$q$  & $(\rho+1+(\beta-1)\alpha)p + (\sigma+1)q-pq$ & -1 & -1 \\
\hline  
$\cdots$ & $\cdots$ & -1 & -1 \\
\hline
$q+p-\sigma-2$ & $(\rho+1+(\beta+\sigma+1-p)\alpha)p+(p-1)q-pq$ & -1 & -1\\
\hline
$q+p-\sigma-1$ & $(\rho+1+(\beta+\sigma-p)\alpha)p$ & 0 & 0\\
\hline
$\cdots$ & $\cdots$ & 0 & 0 \\
\hline
$q+\beta-1$ & $(\rho+1)p+(\beta+\sigma-p)q$ & 0 & 0\\
\hline
$q+\beta$ & $(\rho+1-\alpha)p+(\beta+\sigma+1-p)q$ & 1 & 0\\
\hline
$\cdots$ & $\cdots$ & 1 & 0 \\
\hline
$q+p-1$ & $(\rho+1-(p-\beta)\alpha)p+\sigma q$ & 1 & 0\\
\hline
\end{tabular}
\end{center}
Since $(\rho+1)p+(\sigma+1)q=qp+1$ we can rewrite $w_{+}$ as
$$w_{+}=(\rho+1-(p-\beta)\alpha)p+(\sigma-\beta)q.$$
The condition (\ref{b1})
ensures that $\rho+1-(p-\beta)\alpha \ge 0$. {}From $\beta\le (p-3)/2$ and
$\beta+\sigma\ge p-1$ we infer that $\sigma\ge \beta+2$.
It follows that
$$0\le w_{+}=(p-1)q+1-(p-\beta)\alpha p-\beta q<pq.$$  Thus $f(p-1)=w_{+}$.
The condition (\ref{a1}) ensures that $\rho+1+(\beta-1)\alpha\le q-1$.
Let $0\le r_1<pq$ be the unique integer with $r_1\equiv -{1\over r}({\rm mod~}pq)$. It
is easy to see that $r_1=q-\alpha p$. We have $f(m)\equiv w_{+}+(m-p+1)r_1({\rm mod~}pq)$.
Using these observations one arrives at Table 6A. For $m\le p-\sigma+\beta-2$ we do not care about
whether $\tau_m=0$ or $\tau_m=1$; in either case we find $b_m=0$
and hence $c_m=0$.\\
\indent On noting that $f(m+q)\equiv f(m)+\beta q({\rm mod~}pq)$ (cf. the
proof of Theorem \ref{11}), we easily infer that the $f(m+q)$ are as given
in Table 6B, with the caveat that the entries from $q+p-\sigma-1$ to $q+\beta-1$ do
not arise if $\beta+\sigma=p-1$. Using that $f(q+p-\sigma-2)\le f(p-1)$
and $f(q+\beta)>f(p-1)$ (a consequence of $\alpha$ satisfying (\ref{c1}), 
respectively (\ref{d1}))), we deduce that the
$b_{m+q}$ and $c_{m+q}$ columns are as given in Table 6B.
Finally on applying Lemma \ref{kapel} we infer that
$$a_{pqr}(p-1+rw_{+})=\sum_{m=p-\sigma+\beta-1}^{p-1}1+\sum_{m=q}^{q+p-\sigma-2}1
=(\sigma-\beta+1)+(p-\sigma-1)=p-\beta.$$
This concludes the proof.\qed

\begin{Lem}
\label{abcde2}
Let ${\cal I}$ be the set of real numbers satisfying {\rm (\ref{a1}), 
(\ref{b1}), (\ref{c1}), (\ref{d1})} and suppose the conditions
of Lemma {\rm \ref{main1}} preceding {\rm (\ref{a1})} are satisfied. 
Put
$$\gamma=\min \Big\{{p-\beta^*\over p-\beta},{\beta^*-\beta\over \beta^*} \Big\}.$$
The set ${\cal I}$
is non-empty iff $\beta^*\le 2\beta $. In that case
$${\cal I}=\Big({q(p-1-2\beta)\over p(p-1-\beta)},{q\gamma\over
p}\Big].$$
If ${\cal I}$ is non-empty then it consists of positive reals only.
\end{Lem}
{\it Proof}. Left to the reader. \qed\\

\noindent {\tt Remark}. Note that
$$\gamma=\cases{{p-\beta^*\over p-\beta} & if $p<\beta+{\beta^*\over \beta}(\beta^*-\beta)$;\cr
{\beta^*-\beta\over \beta^*} & otherwise.}$$
\noindent On combining the latter two lemmas one obtains an explicit
counter-example construction in the positive case.
\begin{Thm}
\label{positive}
Suppose that ${\cal B}_{+}(p)$ is non-empty and $\beta \in {\cal B}_{+}(p)$.
Let $q>p$ be a prime satisfying $q\equiv \beta({\rm mod~}p)$ and
$q\ge q_{+}(p)$ with
$$q_{+}(p)={p(p-1-\beta)\over \gamma(p-1-\beta)-p+1+2\beta }.$$ Then the interval
$${\cal I}=\Big[{q(p-1-2\beta)\over p(p-1-\beta)},{q\gamma\over
p}\Big]$$
contains at least one integer $a$.
Let $r>q$ be a prime with
$r(q-pa)\equiv -1({\rm mod~}pq)$,
then
$$a_{pqr}(p-1+[(p-\beta-1)q-(p-\beta)ap]r)=p-\beta>{(p+1)\over 2}$$
is a counter-example to Beiter's conjecture.
\end{Thm}
{\it Proof}. Note that as a function of $q$ the length of ${\cal I}$ is
increasing. If the length of ${\cal I}$ is at least one, then it contains
at least one positive integer. Now $q_{+}(p)$ is obtained on solving
the equation $|{\cal I}|=1$ for $q$. The proof is completed on combining 
Lemma \ref{main1} and Lemma \ref{abcde2}. \qed\\

\noindent Remark. If $1\le \beta\le (p-3)/2$ and $\beta+\beta^*\ge p$, then
$\beta$ satisfies the conditions of Lemma \ref{main1}. If, in addition, $\beta^*\le 2\beta$, that
is if $\beta\in {\cal B}_{+}(p)$, then the conditions of 
both Lemma \ref{main1} and Lemma \ref{abcde2} are satisfied by $\beta$.

\section{The proofs of Theorem \ref{main} and Theorem \ref{main3}}
As is well-known the distribution of inverses modulo $p$ can be studied
by connecting this problem to Kloosterman sums and estimates for those. For
us, the following typical lemma, see e.g. Cobeli \cite[Lemma 4, Chapter 3.2]{Cobeli}, will do.\\
\indent Let ${\cal I}=\{a,a+h, \ldots ,a+(M-1)h \}\subset [1,p]$
and put $$N({\cal I}_1,{\cal I}_2;p)=\#\{(x,y):x\in {\cal I}_1,~y\in {\cal I}_2,~
xy\equiv 1({\rm mod~}p)\},$$
where ${\cal I}_1$ and ${\cal I}_2$ are allowed to have different increments $h$.
\begin{Lem}\label{cobel}
Let $p$ be a prime number. 
We have
$$
\Big| N({\cal I}_1,{\cal I}_2;p)-\frac{|{\cal I}_1|\cdot|{\cal I}_2|}{p}\Big|\leq 
\sqrt{p}~(2+\log p)^2.
$$
\end{Lem}
The set of points $(x,y)$ with $xy\equiv d({\rm mod~}p)$ is called a {\it modular
hyperbola}. For a survey of this area of study see e.g. Shparlinski \cite{S}.\\

\noindent {\it Proof of Theorem} \ref{main3}. Bachman's upper bound $p-\lceil p/4\rceil$ shows
that $M(p)\le 3p/4$.\\ 
\indent Suppose that for $p\ge 29$ and some $0<\epsilon<1/6$ we have
\begin{equation}
\label{starrie}
{p\over 3}(1+\epsilon)\le \beta\le {p\over 3}(1+2\epsilon),~
{2p\over 3}(1-{\epsilon\over 2})\le \beta^*\le {2p\over 3}(1+\epsilon).
\end{equation}
Then one checks that $\beta\in {\cal B}_{+}(p)$. It then follows
by Theorem \ref{main} that $M(p)\ge 2p(1-\epsilon)/3$. It only
remains to show that for every $p$ sufficiently large there
is a $\beta$ satisfying (\ref{starrie}). This follows on invoking
Lemma \ref{cobel} with ${\cal I}_1$ the integers in the interval
$[{p\over 3}(1+\epsilon),{p\over 3}(1+2\epsilon)]$ and ${\cal I}_2$
the integers in the range $[{2p\over 3}(1-{\epsilon\over 2}),{2p\over 3}(1+\epsilon)]$. 
This completes the proof.\qed

\begin{Prop}
\label{duke} Let $\epsilon>0$.
There are infinitely many primes $p$ such that there exist primes $q$ and $r$ so
that $$A(pqr)=\min({p-1\over 2}+a,p-a)\ge ({2\over 3}-\epsilon)p,$$
with $a$ as in Theorem {\rm \ref{flup2}}.
\end{Prop}
{\it Proof}. Suppose that
\begin{equation}
\label{story}
p\ge 29,~0<\epsilon<{1\over 9},~\beta\in {\cal B}_{+}(p),~\beta<{p\over 3}(1+3\epsilon),~\beta+\beta^*=p,
\end{equation} 
then
by Theorem \ref{main} we obtain $a_{pqr_{+}}(n_{+})=p-\beta\ge ({2\over 3}-\epsilon)p$. In this case also
Bachman's upper bound given in Theorem \ref{flup2} gives $A(pqr_{+})=\beta^*=p-\beta$. To see 
this note that $r_{+}^*=p-\beta$ and $q^*=\beta^*$. Since $\beta<\beta^*$ we infer
that $a=p-\beta^*$. By Theorem \ref{flup2} it then follows that
$A(pqr_{+})\le \min({p-1\over 2}+p-\beta^*,\beta^*)\le \beta^*$. Since 
$a_{pqr_{+}}(n_{+})=\beta^*$, it follows that $A(pqr_{+})=\beta^*$ and thus
Bachman's upper bound is assumed.\\
\indent Duke et al. \cite {DFI} proved that if $f$ is a quadratic polynomial with
complex roots, $0\le a<b\le 1$, then
$$\#\{(p,\nu):p\le x,~f(\nu)\equiv 0({\rm mod~}p),~a\le {\nu\over p}<b\}\sim (b-a)\pi(x),$$
where $\pi(x)$ denotes the number of primes $p\le x$. In particular it follows that there are 
asymptotically $\epsilon \pi(x)$ primes
$p$ for which there exist $v$ satisfying $v+v^*=p$ (that is $v^2+1\equiv 0({\rm mod~}p)$),
and $p/3<v<p(1+3\epsilon)/3$. On putting $v=\beta$ we then see that $\beta\in {\cal B}_{+}(p)$ 
and thus it follows that there exist infinitely many primes $p$ for which there is a
$\beta$ satisfying (\ref{story}).\qed\\

\noindent The final result in this section shows that by elementary methods one can easily
prove that $M(p)\ge ({2\over 3}-\epsilon)p$ for infinitely many primes $p$. 
\begin{Prop}
\label{yves}
Let $\epsilon>0$, $e \geq 1$ be an integer, $N = 2^{2e+1}$ and $p$
a prime satisfying $p \equiv N - 9~({\rm mod~}3N)$. If $p \ge {N^2\over 2} - 9$,
then $$M(p)\geq {(2N-1)p - 9\over 3N} > \Big({2\over 3} - {N\over 3(N^2 - 18)}\Big)p.$$
If  $p \ge {N^2\over 2} - 9$ and $N>{1\over 3\epsilon}+3$, then
$M(p)>({2\over 3}-\epsilon)p$.
\end{Prop}
{\it Proof}. Consider $\beta =((N+1)p + 9)/(3N)$. Then $\beta^* =(2p+N)/3$.
One checks that if $p \ge N^2/2 - 9$, then $\beta^*/2 \le \beta$ and
finally that $\beta \in {\cal B}_{+}(p)$. Then invoke Theorem \ref{main}. \qed\\

\noindent {\tt Examples}.\\
If $p \equiv 23({\rm mod~}24)$ and $p > 23$ then $M(p) \geq (5p - 3)/8 > 0.608p$.\\
If $p \equiv 23({\rm mod~}96)$ and $p > 503$ then $M(p) \geq (21p - 3)/32 > 0.656p$.\\
If $p \equiv 119({\rm mod~}384)$ and $p > 8183$ then $M(p) \geq (85p - 3)/128 > 0.664p$.\\
If $p \equiv 503({\rm mod~}1536)$ and $p > 131063$ then $M(p) \geq (341p - 3)/512 > 0.666p$.\\

\noindent {\it Proof of Theorem} \ref{main}. Follows on combining Theorem
\ref{negative} and Theorem \ref{positive}. Theorem \ref{positive} together
with Theorem \ref{flup2} yields that 
$A(pqr_{+})=p-\beta$ if $\beta\in {\cal B}_{+}(p)$ and $\beta+\beta^*=p$
(cf. the proof of Proposition \ref{duke}). \qed

\section{Reciprocal cyclotomic polynomials}
We point out that for the so called reciprocal cyclotomic polynomials the
analogue of Beiter's conjecture is known.
\noindent Let 
$${1\over \Phi_n(x)}=\sum_{k=0}^{\infty}c_n(k)x^k$$
be the Taylor series of $1/\Phi_n(x)$ around $x=0$. 
The coefficients turn out to be periodic with period dividing $n$.
Moree \cite{Mor} established
the following result concerning the height, $H(n)$, of 
$1/\Phi_n(x)$ (thus $\max_{k\ge 0}|c_n(k|=H(n)$).
\begin{Thm}
Let $p<q<r$ be odd primes. Then
$H(pqr)=p-1$ iff 
$$q\equiv r\equiv \pm 1({\rm mod~}p){\rm ~and~}r<{(p-1)\over (p-2)}(q-1).$$
In the remaining cases $H(pqr)<p-1$.
\end{Thm} 
This result in
combination with Dirichlet's theorem on arithmetic progressions shows that
for every odd prime $p$ there are infinitely many pairs $(q,r)$ such that
$H(pqr)=p-1$.\\
\indent Let $m$ be an arbitrary natural number. In \cite{JLM} simple properties of reciprocal cyclotomic
polynomials are used to show that $\{a_{mn}(k)~|~n\ge 1,~k\ge 0\}=\mathbb Z$
and, likewise, $\{c_{mn}(k)~|~n\ge 1,~k\ge 0\}=\mathbb Z$.\\

\noindent {\tt Acknowledgement}. We thank C. Cobeli, 
M.Z. Garaev and I. Shparlinski for helpful information concerning the
distribution of inverses modulo $p$ and the two referees for their
careful proofreading.
N. Kaplan pointed out to
us that in the summer of 2007 Tiankai 
Liu (Harvard) wrote a program that computed some counter-examples. We acknowledge
that we were not the first to find counter-examples.
In Moree \cite{Mor} it is shown that the analogue of Beiter's conjecture is
false for reciprocal cyclotomic ternary polynomials. This gave us the idea (not
being aware of Liu's work)
to do a thorough numerical check on the original Beiter conjecture, leading
to our first counter-example on Sept. 9, 2007.

\medskip\noindent {\footnotesize 12 bis rue Perrey,\\
31400 Toulouse, France.\\
e-mail: {\tt galloty@orange.fr}}\\

\medskip\noindent {\footnotesize Max-Planck-Institut f\"ur Mathematik,\\
Vivatsgasse 7, D-53111 Bonn, Germany.\\
e-mail: {\tt moree@mpim-bonn.mpg.de}}
\vskip 5mm

\end{document}